\documentclass{siamart171218}

\usepackage{amsfonts}
\usepackage{graphicx}
\usepackage{hyperref}
\usepackage{amsmath,amssymb,setspace}
\usepackage{algorithm}
\usepackage{graphics}
\usepackage{booktabs}
\usepackage{algorithmic}
\usepackage{bm}
\newcommand{\bs}{\boldsymbol{s}}
\newcommand{\bH}{\boldsymbol{H}}
\newcommand{\bd}{\boldsymbol{d}}
\newcommand{\bL}{\boldsymbol{L}}
\newcommand{\bu}{\boldsymbol{u}}
\newcommand{\bP}{\boldsymbol{P}}
\newcommand{\bp}{\boldsymbol{p}}
\newcommand{\bV}{\boldsymbol{V}}
\newcommand{\bv}{\boldsymbol{v}}
\newcommand{\bM}{\boldsymbol{M}}
\newcommand{\bU}{\boldsymbol{U}}
\newcommand{\bT}{\boldsymbol{T}}
\newcommand{\bq}{\boldsymbol{q}}
\newcommand{\be}{\boldsymbol{e}}

\headers{An Example Article}{J. Jiang, F. Sanogo, and C. Navasca}

\title{Low-CP-rank Tensor Completion via Practical Regularization  \thanks{Submitted to the editors January 15, 2021.
\funding{This work was funded by National Science Foundation under Grant No. DMS-1439786.}}}

\author{Jiahua Jiang\thanks{School of Information Science and Technology, ShanghaiTech University, Shanghai, China 
  (\email{jiangjh@shanghaitech.edu.cn)}.}
\and Fatoumata Sanogo\thanks{Department of Mathematics, University of Alabama at Birmingham, Birmingham, AL 35294
  (\email{sanogof1@uab.edu)}.}
\and Carmeliza Navasca\thanks{Department of  Mathematics, University of Alabama at Birmingham, Birmingham, AL 35294
  (\email{cnavasca@uab.edu)}.}}

\begin{document}

\maketitle

% REQUIRED
\begin{abstract}
Dimension reduction techniques are often used when the high-dimensional tensor has relatively low intrinsic rank compared to the ambient dimension of the tensor. The CANDECOMP/PARAFAC (CP) tensor completion is a widely used approach to find a low-rank approximation for a given tensor.  In the tensor model \cite{8645405}, an $\ell_1$ regularized optimization problem was formulated with an appropriate choice of the regularization parameter. The choice of the regularization parameter is important in the approximation accuracy. However, the emergence of the large amount of data poses onerous computational burden for computing the regularization parameter via classical approaches \cite{gazzola2020krylov} such as the weighted generalized cross validation (WGCV) \cite{chung2008weighted}, the unbiased predictive risk estimator \cite{stein1981estimation, vogel2002computational}, and the discrepancy principle \cite{morozov1966solution}. In order to improve the efficiency of choosing the regularization parameter and leverage the accuracy of the CP tensor, we propose a new algorithm for tensor completion by embedding the flexible hybrid method \cite{gazzolaflexible} into the framework of the CP tensor. The main benefits of this method include incorporating regularization automatically and efficiently, improved reconstruction and algorithmic robustness.  Numerical examples from image reconstruction and model order reduction demonstrate the performance of the propose algorithm. 
\end{abstract}

% REQUIRED
%\begin{AMS}
 % 68Q25, 68R10, 68U05
%\end{AMS}

\section{Introduction}
Tensor computations have become prevalent in across many fields in mathematics \cite{2006math7647D,105555270146}, computer science \cite{2017arXiv171104887G,8408516,7931588}, engineering \cite{chadwick1995} and data science \cite{Hou2017TensorbasedRM,1010071142799521}. In particular, tensor methods are now ubiquitous in areas of numerical linear algebra \cite{BECKMANN2005294,Beylkin05algorithmsfor}, imaging science \cite{8319500,8451531} and applied algebraic geometry \cite{KoBaKe05}. New tensor based methods are currently being developed in scientific computing of complex problems \cite{BOELENS2020109744}. %Cite many papers here by D. Venturi, B. Khoromskij, and etc..

The tensor rank problem is crucial in reconstructing a given tensor $\mathcal{T}$ through a sum of rank-one tensor product, $\sum_{r}^{R} \alpha_r\mathbf{a}_{r} \otimes \mathbf{b}_{r} \otimes \mathbf{c}_{r}$. In practice, the column vectors, $\mathbf{a}_{r}$,  $\mathbf{b}_{r}$ and  $\mathbf{c}_{r}$, are concatenated into what we call factor matrices, $\mathbf{A}$, $\mathbf{B}$ and $\mathbf{C}$. The minimum number of summands $R$ is the rank of the tensor. The elements of the vector $\sigma$ of size $R$
are the scalings of the rank-one tensors. This tensor factorization is the well-known canonical polyadic or CANDECOMP/PARAFAC (CP) decomposition. Optimization strategies like the alternating least squares and gradient descent algorithms are some numerical methods are the standard methods given an input of the tensor rank $R$. 

A tensor model \cite{8645405} which incorporates tensor rank approximation is the following:
\begin{eqnarray*}\label{opt1}
\min_{A,B,C,\sigma} \left \Vert \mathcal{T} - \mathcal{S} \right \Vert_F + \lambda \Vert \alpha \Vert_{\ell_1}\\
\end{eqnarray*}
where $\lambda$ is a regularization parameter and $\mathcal{S} = \sum_{r}^{U} \alpha_r \mathbf{a}_{r} \otimes \mathbf{b}_{r} \otimes \mathbf{c}_{r}$ with an upper bound tensor rank $R$. This sparse optimization problem is solved iteratively and the vector $\alpha$ reveals an approximated tensor rank $\Vert \gamma \Vert_{\ell_0}=R$ with $R << U$. The drawback of this model is its dependent on the regularization parameter $\lambda$. In the work \cite{WangNa2}, the choice of the regularization parameter is tied to two intrinsic parameters: variance of noise and incoherence of the given tensor data $\mathcal{T}$. One has to initialize $\lambda$ from a bound based on these two parameters from the tensor data with an upper bound rank $R$. In practice, a priori estimates of the variance and incoherence parameters are needed based on a CP decomposition of the given data with an initial tensor rank guess $R$. Only then $\lambda$ can be chose accordingly from an estimated bound. The advantage of this approach provides the theoretical bounds for $\lambda$. However, it is not practical enough for real data implementation; the choice of $\lambda$ is only as good as estimated intrinsic parameters. Moreover, the accuracy level is only around $10^{-2}$. In this paper, we show a more adaptive, practical and methodical way for calculating the regularization parameter $\lambda$ using the flexible hybrid method. The flexible hybrid method allows more efficiency than classical approaches \cite{gazzola2020krylov}, such as, the weighted generalized cross validation (WGCV) \cite{chung2008weighted}, the unbiased predictive risk estimator \cite{stein1981estimation, vogel2002computational}, and the discrepancy principle \cite{morozov1966solution} in problems with large amount of data. We have shown that new iterative method gives more accurate results in tensor completion and model order reduction problems.

We have two application areas in this paper: tensor completion in image restoration and model order reduction. Matrix and tensor completion techniques provide major tools in recommender systems in computer science and in general in data science; it is about filling in missing entries from the partially observed entries of the matrix or tensor. The success of matrix completion methods are attributed to sparse optimization methods in compressed sensing \cite{candes2009power}. These methods have been generalized iteratively to tensor completion problem \cite{6138863} via the matricized tucker models \cite{WALCZAK200115,ANDERSSON199893} where missing entries are predicted through the trace norm optimization. In fact, the tensor completion problem dates back as early as in 2000. Bro \cite{BRO1997149} had one of the earliest work on demonstrating two ways to handle missing data using CP. The first way is to alternatively estimate the model parameters while imputing the missing data. Another approach called Missing-Skipping skips the missing value and builds up the model based only on the observed part via a weighted least squares formulation in the CP format \cite{AcDuKo11}. Our proposed tensor completion gives more accurate results in capturing more features in color images through low rank construction via our model and numerical technique.

Furthermore, to alleviate the computational effort for many-query computations and repeated output evaluations for different values of some inputs of interest, besides classical model order reduction approaches \cite{ohlberger2016approximation, nouy2015low} such as Reduced Basis Methods \cite{binev2011convergence, rozza2007reduced} and Proper Orthogonal Decomposition (POD), much recent effort in tensor-based model reduction such as Randomized CP tensor decomposition \cite{erichson2020randomized} and tensor POD \cite{zhang2017design}, has been rewarded with many
promising developments. Compared with the classical model order reduction approaches, tensor-based model reduction algorithms allow us to achieve significant computational savings, especially for expensive high fidelity numerical solvers. 
% {\textcolor{red}{Jiaha: can you add two-three sentences on model order reduction related the sections on this paper.}}

\subsection{Organization}

The paper is organized as follows. In Section \ref{sec:main}, we provide some tensor backgrounds and basic tensor decomposition.
Then, Section \ref{sec:3} deals with deriving the iterative equation by through matricizations and alternating block optimization incorporating the proximal gradient formulation and the flexible hybrid method for the automatic selection of the regularization parameter. Experimental
results are in Section \ref{sec:experiments}, and the conclusions follow in
Section \ref{sec:conclusions}.

\section{Preliminaries}
\label{sec:main}
We denote a vector by a bold lower-case letter $\mathbf{a}$.
The bold upper-case letter $\mathbf{A}$ represents a matrix
and the symbol of tensor is a calligraphic letter $\mathcal{A}$.
Throughout this paper, we focus on third-order tensors
$\mathcal{A}=(a_{ijk})\in\mathbb{R}^{I\times J\times K}$ of three indices $1\leq i\leq I,1\leq j\leq J$ and $1\leq k\leq K$,
but all are applicable to tensors of arbitrary order greater or equal to three.

A third-order tensor $\mathcal{A}$ has column, row and tube fibers,
which are defined by fixing every index but one and
denoted by $\mathbf{a}_{:jk}$, $\mathbf{a}_{i:k}$ and $\mathbf{a}_{ij:}$ respectively.
Correspondingly, we can obtain three kinds $\mathbf{A}_{(1)},\mathbf{A}_{(2)}$ and $\mathbf{A}_{(3)}$ of matricization
of $\mathcal{A}$ according to respectively arranging the column, row, and tube fibers to be columns of matrices.
We can also consider the vectorization for $\mathcal{A}$ to
obtain a row vector $\mathbf{a}$ such the elements of $\mathcal{A}$ are arranged
according to $k$ varying faster than $j$ and $j$ varying faster than $i$, i.e.,
$\mathbf{a}=(a_{111},\cdots,a_{11K},a_{121},\cdots,a_{12K},\cdots,a_{1J1},\cdots,a_{1JK},\cdots)$.

Henceforth, the outer product of a rank-one third order tensor is denoted as $\mathbf{a}\circ\mathbf{b}\circ\mathbf{c}\in\mathbb{R}^{I\times J\times K}$ of three nonzero vectors $\mathbf{a}, \mathbf{b}$ and $\mathbf{c}$ is a rank-one tensor with elements $a_ib_jc_k$ for all the indices.
A canonical polyadic decomposition of $\mathcal{T}\in\mathbb{R}^{I\times J\times K}$ expresses $\mathcal{T}$ as a sum of rank-one outer products:
\begin{equation}\label{cpd}
\mathcal{T}=\sum_{r=1}^{R} \mathbf{a}_r\circ\mathbf{b}_r\circ\mathbf{c}_r
\end{equation}
where $\mathbf{a}_r\in\mathbb{R}^I,\mathbf{b}_r\in\mathbb{R}^J,\mathbf{c}_r\in\mathbb{R}^K$ for $1\leq r\leq R$.
Every outer product $\mathbf{a}_r\circ\mathbf{b}_r\circ\mathbf{c}_r$ is called as a rank-one component
and the integer $R$ is the number of rank-one components in tensor $\mathcal{A}$.
The minimal number $R$ such that the decomposition (\ref{cpd}) holds is the rank of tensor $\mathcal{A}$, which is denoted by $\mbox{rank}(\mathcal{A})$.
For any tensor $\mathcal{A}\in\mathbb{R}^{I\times J\times K}$,
$\mbox{rank}(\mathcal{A})$ has an upper bound $\min\{IJ,JK,IK\}$ \cite{kruskal}.

In this paper, we consider CP decomposition in the following form
\begin{equation}\label{cpd2}
\mathcal{T}=\sum_{r=1}^{R} \alpha_r\mathbf{a}_r\circ\mathbf{b}_r\circ\mathbf{c}_r
\end{equation}
where $\alpha_r\in\mathbb{R}$ is a rescaling coefficient of rank-one tensor $\mathbf{a}_r\circ\mathbf{b}_r\circ\mathbf{c}_r$ for $r=1,\cdots,R$.
For convenience, we let ${\alpha}=(\alpha_1,\cdots,\alpha_R)\in\mathbb{R}^R$
and denote $[{\alpha};\mathbf{A},\mathbf{B},\mathbf{C}]_R = \sum_{r=1}^{R} \alpha_r\mathbf{a}_r\circ\mathbf{b}_r\circ\mathbf{c}_r$
in (\ref{cpd2}) where $\mathbf{A}=(\mathbf{a}_1,\cdots,\mathbf{a}_R)\in\mathbb{R}^{I\times R},\mathbf{B}=(\mathbf{b}_1,\cdots,\mathbf{b}_R)\in\mathbb{R}^{J\times R}$ and $\mathbf{C}=(\mathbf{c}_1,\cdots,\mathbf{c}_R)\in\mathbb{R}^{K\times R}$
are called the factor matrices of tensor $\mathcal{A}$.

In most iterative techniques for decompositions, tensor \emph{matricizations} are required to transform the tensor equations into matrix equations. Here we describe a standard approach for a matricizing of a tensor. The Khatri-Rao product \cite{Smilde2004MultiwayAW} of two matrices $X\in\mathbb{R}^{I\times R}$ and $Y \in\mathbb{R}^{J\times R}$
is defined as
$${X\odot Y}=({x}_1\otimes {y}_1,\cdots, {x}_R \otimes {y}_R)\in {R}^{IJ\times R},$$
where the symbol ``${\otimes}$" denotes the Kronecker product:
$${x\otimes y}=(x_1y_1,\cdots,x_1y_J,\cdots,x_Iy_1,\cdots,x_Iy_J)^T.$$
Using the Khatri-Rao product, the CP model \cref{cpd2} can be written in three equivalent matrix equations:
\begin{subequations}
\begin{align}
  \mathbf{T}_{(1)}
  &= \mathbf{AD}(\mathbf{C\odot B})^T, \label{kathri1} \\
  \mathbf{T}_{(2)}
  &=\mathbf{BD}(\mathbf{C\odot A})^T,  \label{kathri2} \\
  & \mbox{and} \nonumber \\
  \mathbf{T}_{(3)}
  &=\mathbf{CD}(\mathbf{B\odot A})^T \label{kathri3}
\end{align}
\end{subequations}
where the matrix $\mathbf{D}$ is diagonal with elements of ${\alpha}$.  To achieve CP decomposition of given tensor $\mathcal{T}$ with a known tensor rank $R$ and an assumption that $\mathbf{D}=\mathbf{I}$, the matrix equations \cref{kathri1}-\cref{kathri3} are formulated into linear least-squares subproblems to solve iteratively for $\mathbf{A}$, $\mathbf{B}$ and $\mathbf{C}$, respectively. Here are the linear least-squares subproblems:
\begin{subequations}
\begin{align}
  &\min\limits_{\mathbf{A}}\Vert \mathbf{T}_{(1)} - \mathbf{AD}(\mathbf{C\odot B})^T \Vert_F^2 , \label{kathri1b} \\
  &\min\limits_{\mathbf{B}} \Vert \mathbf{T}_{(2)} - \mathbf{BD}(\mathbf{C\odot A})^T \Vert_F^2,  \label{kathri2b} \\
   & \mbox{and} \nonumber \\
  & \min\limits_{\mathbf{C}}  \Vert \mathbf{T}_{(3)} - \mathbf{CD}(\mathbf{B\odot A})^T \Vert_F^2. \label{kathri3b}
\end{align}
\end{subequations}
This technique is the well known Alternating Least-Squares (ALS) \cite{BRO1997149,li2011convergence}. Typically, a normalization constraint on factor matrices such that each column is normalized to length one \cite{doi:10.1137/110843587,Acar_2011} is required for convergence, which we denote by $\mathbf{N(A,B,C)}=1$.
%Add convergence discussion, Uschmajew, paper with Na Li.

\section{Iterative Equations for Tensor Completion}
\label{sec:3}

We will describe our low rank tensor model of a given tensor in a CP format with an approximated tensor rank for tensor completion.
Our goal is to fill in the missing entries from a given tensor $\mathcal{T}$ with the partially observed entries by reconstructing a completed low rank tensor $\mathcal{S}$. To do so, we formulate a sparse optimization problem \cite{8645405} for recovering CP decomposition from tensor $\mathcal{T} \in \mathbb{R}^{I \times J \times K}$ with partially observed entries on the index set $\Omega$: 
\begin{subequations}
\begin{align}
  &\min_{A,B,C,\sigma} \left \Vert \mathcal{T} - \mathcal{S} \right \Vert_F + \lambda \Vert \alpha \Vert_{\ell_1} \label{prob}\\
  &\mbox{subject~to~}\mathcal{S}(\Omega)=\mathcal{T}(\Omega) \nonumber
\end{align}
\end{subequations}
where $\lambda$ is a constant regularization parameter and  $\mathcal{S}= \sum_r \alpha    
\mathbf{a}_r\circ\mathbf{b}_r\circ\mathbf{c}_r$.

We will now derive the iterative equations for $A, B, C$ and $\sigma$. The equations are typically associated with Iterative Soft Thresholding Algorithm (ISTA) \cite{Beck2009AFI} whose derivation is based on the Majorization-Minimization (MM) \cite {4358845} method. ISTA (Iterative Soft-Thresholding Algorithm) is a combination of the Landweber algorithm and soft-thresholding (so it is also called the thresholded-Landweber algorithm).% This algorithm can be traced back to the proximal forward-backward iterative scheme introduced by Bruck (1977) and Passty (1979) within the general framework of splitting methods.

%For the Majorization-Minimization (MM) we pick a vector $x_k$, a `guess' for the minimum of $J(x)= \left \Vert \mathcal{T} - \mathcal{S} \right \Vert_F + \lambda \Vert \sigma \Vert_{\ell_1}$.Based on $x_k$, we want to  find $x_{k+1}$ such that $J(x_{k+1}) < J(x_k)$. The MM approach is first to choose a new function (G(x))which can easily be minimize such that $G(x)\geq J(x)$ for all x and minimize G(x) to get $x_{k+1}$. G(x) should equal J(x) at $x_k$.
%The tensor completion problem  (\ref{prob}) is now a minimization of a sum of a smooth term and a nonsmooth term

Suppose we have a minimization problem:
\begin{equation}\label{pg1}
\min\limits_{{x}} f(x). 
\end{equation}
By using the proximal operators formulation (see Appendix) and the MM approach, we first find an upper bound for $f(x)$:
\begin{equation*}
 f(x) \leq f(y) + \langle \nabla_x f(y) , x-y \rangle + \gamma \Vert x - y \Vert^2_2 
\end{equation*}
Let $g(x,y)=f(y) + \langle \nabla_x f(y) , x-y \rangle + \gamma \Vert x - y \Vert^2_2$. Note that $f(x) \leq g(x,y)$ for all $x$ and $f(x) = g(x,y)$ when $y=x$.
Thus, we can reformulate \ref{pg1} as
\begin{equation}\label{pg2}
\min\limits_{{x}} f(y) + \langle \nabla_x f(y) , x-y \rangle + \gamma \Vert x - y \Vert^2_2.
\end{equation}
Since this is a minimization over $x$, then \ref{pg2} is equivalent to
\begin{equation}\label{pg3}
\min\limits_{{x}} \langle \nabla_x f(y) , x-y \rangle + \gamma \Vert x - y \Vert^2_2.
\end{equation}
By gathering the terms with respect to $x$, the objective function in \ref{pg3} can be expressed as
\begin{equation*}
     \gamma (-2b^Tx + x^Tx) + c,
\end{equation*}
where $c=\gamma (y)^Ty - \nabla f(y)^Ty$ and $b = y - \frac{1}{2 \gamma} \nabla_x f(y)$.
Since $b^Tb - 2b^Tx + x^Tx= \Vert b -x \Vert^2_2$, we have a new formulation:
\begin{equation}\label{pg4}
\min\limits_{{x}} \gamma \Vert x - b \Vert_2^2.
\end{equation}
Now from the least-squares problems (\ref{kathri1b}-\ref{kathri3b}) and using proximal gradient formulation, we have the following new formulations:
\begin{equation*}
\begin{aligned}
\mathbf{A}^*=
& \arg\min\limits_{\mathbf{A}} \{\langle\mathbf{A}-\mathbf{A}^n,\nabla_\mathbf{A} f(\mathbf{A}^n,\mathbf{B}^n,\mathbf{C}^n,{\alpha}^n)\rangle+\frac{sd_n}{2}\|\mathbf{A}-\mathbf{A}^n\|_F^2\} &\\
& \quad\mbox{s.t.}\ \|\mathbf{a}_i\|=1, i=1,\cdots,R, &
\end{aligned}
\end{equation*}
\begin{equation*}
\begin{aligned}
\mathbf{B}^*=
& \arg\min\limits_{\mathbf{B}} \{\langle\mathbf{B}-\mathbf{B}^n,\nabla_\mathbf{B} f(\mathbf{A}^{n+1},\mathbf{B}^n,\mathbf{C}^n,{\alpha}^n)\rangle+\frac{se_n}{2}\|\mathbf{B}-\mathbf{B}^n\|_F^2\}&\\
& \quad\mbox{s.t.}\ \|\mathbf{b}_i\|=1, i=1,\cdots,R, &
\end{aligned}
\end{equation*}
and
\begin{equation*}
\begin{aligned}
\mathbf{C}^*=
& \arg\min\limits_{\mathbf{C}} \{\langle\mathbf{C}-\mathbf{C}^n,\nabla_\mathbf{C} f(\mathbf{A}^{n+1},\mathbf{B}^{n+1},\mathbf{C}^n,{\alpha}^n)\rangle+\frac{sf_n}{2}\|\mathbf{C}-\mathbf{C}^n\|_F^2\}&\\
& \quad\mbox{s.t.}\ \|\mathbf{c}_i\|=1, i=1,\cdots,R, &
\end{aligned}
\end{equation*}
The gradients of $f(\mathbf{A,B,C},\alpha)=\frac{1}{2} \Vert \mathcal{T} - \mathbf{a}_r\circ\mathbf{b}_r\circ\mathbf{c}_r \Vert_F^2$ on $\mathbf{A,B,C,\sigma}$ are the following in terms of the Khatri-Rao product via matricizations:
\begin{subequations}
\begin{align}
    \nabla_\mathbf{A} f(\mathbf{A,B,C,}{\alpha})
    &=(\mathbf{AD}(\mathbf{C\odot B}) - \mathbf{T}_{(1)}) (\mathbf{C\odot B})^T \mathbf{D}^T ,\\
    \nabla_\mathbf{B} f(\mathbf{A,B,C,}{\alpha})
    &=(\mathbf{BD}(\mathbf{C\odot A})-\mathbf{T}_{(2)}) (\mathbf{C\odot A})^T\mathbf{D}^T,\\
    & \mbox{and} \nonumber \\
    \nabla_\mathbf{C} f(\mathbf{A,B,C,}{\alpha}) 
    &=(\mathbf{CD}(\mathbf{B\odot A})-\mathbf{T}_{(3)}) (\mathbf{B\odot A})^T\mathbf{D}^T.
\end{align}
\end{subequations}

Based from the calculations (\ref{pg2}-\ref{pg4}), we obtain the following iterative formula for $A$:
\begin{equation*}
\arg\min\limits_{\mathbf{X}} \{\|\mathbf{X}-\mathbf{D}^n\|_F^2\} \quad\mbox{s.t.}\ \|\mathbf{x}_i\|=1, i=1,\cdots,R.
\end{equation*}
where $\mathbf{D}^n=\mathbf{X}^n-\frac{1}{sd_n}\nabla_\mathbf{X} f(\mathbf{X}^n,\mathbf{Y}^n,\mathbf{Z}^n,{\alpha}^n)$.
We can break it further component-wise: 
\begin{equation*}\label{updatex}
\mathbf{x}_i^{n+1}=\mathbf{d}_i^n/\|\mathbf{d}_i^n\|, i=1,\cdots,R,
\end{equation*}
where $\mathbf{x}_i^{n+1}$ and $\mathbf{d}_i^n$ are the $i$-th columns of $\mathbf{X}^{n+1}$ and $\mathbf{D}^n$.

Similarly, the update of $\mathbf{Y}$ is 
\begin{equation*}
\arg\min\limits_{\mathbf{Y}} \{\|\mathbf{Y}-\mathbf{E}^n\|_F^2\} \quad\mbox{s.t.}\ \|\mathbf{y}_i\|=1, i=1,\cdots,R.
\end{equation*}
where $\mathbf{E}^n=\mathbf{Y}^n-\frac{1}{sd_n}\nabla_\mathbf{Y} f(\mathbf{X}^n,\mathbf{Y}^n,\mathbf{Z}^n,{\alpha}^n)$.
Column-wise, we have
\begin{equation*}\label{updatey}
\mathbf{y}_i^{n+1}=\mathbf{e}_i^n/\|\mathbf{e}_i^n\|, i=1,\cdots,R,
\end{equation*}
where $\mathbf{x}_i^{n+1}$ and $\mathbf{e}_i^n$ are the $i$-th columns of $\mathbf{Y}^{n+1}$ and $\mathbf{E}^n$.

Furthermore, the update of $\mathbf{Z}$ is 
\begin{equation*}
\arg\min\limits_{\mathbf{Z}} \{\|\mathbf{Z}-\mathbf{F}^n\|_F^2\} \quad\mbox{s.t.}\ \|\mathbf{z}_i\|=1, i=1,\cdots,R.
\end{equation*}
where $\mathbf{F}^n=\mathbf{Z}^n-\frac{1}{sd_n}\nabla_\mathbf{Z} f(\mathbf{X}^n,\mathbf{Y}^n,\mathbf{Z}^n,{\alpha}^n)$.
Also, we update vector-wise: 
\begin{equation*}\label{updatez}
\mathbf{z}_i^{n+1}=\mathbf{f}_i^n/\|\mathbf{f}_i^n\|, i=1,\cdots,R,
\end{equation*}
where $\mathbf{z}_i^{n+1}$ and $\mathbf{f}_i^n$ are the $i$-th columns of $\mathbf{Z}^{n+1}$ and $\mathbf{F}^n$.

\subsection{Iterative equation for $\alpha$}

Using the vectorization of tensors in Section $2$,
we can vectorize every rank-one tensor of outer product $\mathbf{a}_r\circ\mathbf{b}_r\circ\mathbf{c}_r$ into a row vector $\mathbf{q}_r$ for $1\leq r\leq R$. We denote a matrix consisting of all $\mathbf{q}_r$ for $1\leq r\leq R$ by
\begin{equation}\label{update2}
\mathbf{Q}=(\mathbf{q}_1^T,\cdots,\mathbf{q}_R^T)^T.
\end{equation}
Thus the function $\frac{1}{2} \Vert \mathcal{T} - \mathbf{a}_r\circ\mathbf{b}_r\circ\mathbf{c}_r \Vert_F^2 $ can be also written as
$\frac{1}{2}\|\mathbf{t}-{\alpha Q}\|_F^2$
where $\mathbf{t}$ is a vectorization for tensor $\mathcal{T}$. Also, the gradient of $f(\bullet)$ on $\mathbf{A,B,C,\alpha}$ is the following in terms of the Khatri-Rao product via matricizations:
\begin{equation*}
 \nabla_{\alpha} f(\mathbf{A,B,C,}{\sigma})=({\alpha}\mathbf{Q}-\mathbf{t})\mathbf{Q}^T.
\end{equation*}
Then, the minimization problem for $\alpha$ is
\begin{equation}\label{alphahOPT}
\min \limits_{\alpha} \frac{1}{2} \|\mathbf{t}-{\alpha \mathbf{Q}}\|_F^2 + \lambda \Vert \alpha \Vert_1.
\end{equation}
Efficiently and appropriately choosing the regularization parameter $\lambda$ plays a crucial role in solving (\ref{alphahOPT}). In the papers \cite{WangNa1,WangNa2}, the proximal operators formulation (see Appendix) and the MM approach are used to solve $\alpha$ iteratively via 
\begin{equation*} 
\alpha^{n+1}=
\arg\min\limits_{{\alpha}} \{\langle{\alpha}-{\alpha}^n,\nabla_{{\alpha}} f(\mathcal{C}^{n+1},\mathbf{A}^{n+1},\mathbf{B}^{n+1},\mathbf{C}^{n+1},{\alpha}^n)\rangle+\frac{s\eta_n}{2}\|{\alpha}-{\alpha}^n\|^2+\lambda\|{\alpha}\|_1\}.
\end{equation*} which is equivalent to the following:
\begin{equation}\label{alpha}
\alpha^{n+1}=\arg\min\limits_{{\alpha}} \frac{1}{2}\|{\alpha}-{\alpha}^n+\frac{1}{s\eta_n}\nabla_{{\alpha}} f(\mathcal{C}^{n+1},\mathbf{A}^{n+1},\mathbf{B}^{n+1},\mathbf{C}^{n+1},{\alpha}^n)\|^2
+\frac{\lambda}{s\eta_n}\|{\alpha}\|_1.
\end{equation}
However, we found that the accuracy of these methods heavily depends on the choice of the initial value of $\alpha$, which reduces the robustness of the whole algorithm, in particular for practical problems. To address this problem, we embed the flexible hybrid method introduced in the following section into the CP completion framework.  
% In the following section, we discuss our approach for solving (\ref{alphahOPT}).

%\begin{equation*}
%\bm{\alpha}^{n+1}=\mathcal{S}_{\frac{\lambda}{s\eta_n}}(\bm{\beta}^{n+1})
%\end{equation*}

% \subsection{Flexible Hybrid Methods for $\ell_1$ regularization}
%\vspace{.2in}

\subsection{The Flexible Hybrid Method for $\ell_1$ Regularization}
The iteratively reweighted norm (IRN) methods \cite{gorodnitsky1992new,rodriguez2008efficient} are typical strategies for solving the $\ell_p-$ regularization inverse problem. However, these methods assume that an appropriate value of the regularization parameter is known in advance, which is hard oftentimes. Therefore, there have been some recent works \cite{giryes2011projected, gazzolaflexible} on selecting regularization parameters for $\ell_p$. In this work, we focus on employing the flexible  method based on Golub-Kahan process \cite{gazzolaflexible} to solve the $\ell_1-$regularized problem,
\begin{equation}
\label{problem:l1}
    \underset{\bs}{\min}||\bH\bs - \bd||^2_2 + \lambda||\bs||_1,
\end{equation}
where $\bd \in \mathbb{R}^m$ is the observed data, $\bH \in \mathbb{R}^{m \times n}$ models the forward process, $\bs \in \mathbb{R}^n$ is the approximation of the desired solution. The first step is to break the $\ell_1-$regularized problem (\ref{problem:l1})  into a sequence of $\ell_2-$norm problems,
\begin{equation}
\label{problem:l2}
    \underset{\bs}{\min}||\bH\bs - \bd||^2_2 + \lambda||\bL(\bs)\bs||^2_2,
\end{equation}
where 
\begin{equation}
    \bL(\bs) = \text{diag}\Big(\big(1/\sqrt{f_{\tau}([|\bs|]_i)}\big)_{i=1,\dots,n} \Big),
\end{equation}
and $f_{\tau}([|\bs|]_i) = \begin{cases}[|\bs|]_i & [|\bs|]_i \geq \tau_1 \\ \tau_2 & [|\bs|]_i < \tau_1\end{cases}$. Here $0 < \tau_2 < \tau_1$ are small thresholds enforcing some additional sparsity in $f_{\tau}([|\bs|]_i)$. Since directly solving (\ref{problem:l2}) is not possible in real problems since the true $\bs$ is not available. To avoid nonlinearities and follow the common practice of iterative methods, $\bL(\bs)$ can be approximated by  $\bL(\bs_k)$, where $\bs_k$ is the numerical solution at the $(k-1)-$th iteration that can be treated as an approximation of the solution at $k-$th iteration. Since directly choosing regularization parameters for large problems is quite costly, the flexible hybrid approaches based on the flexible Golub-Kahan process\cite{gazzolaflexible} has been developed to solve the following variable-preconditioned Tikhonov problem, 
\begin{equation}
    \label{problem:l3}
    \underset{\bs}{\min}||\bH\bs - \bd||^2_2 + \lambda||\bL_k\bs||^2_2,
\end{equation}
which is equivalent to 
\begin{equation}
    \underset{\bs}{\min}||\bH\bL^{-1}_k\widehat{\bs} - \bd||^2_2 + \lambda||\widehat{\bs}||^2_2,
\end{equation}
where $\widehat{\bs} = \bL_k\bs$, and $\bL_k = \bL(\bs_k)$ may change at each iteration. To be able to incorporate the changing preconditioner, the flexible Golub-Kahan process (FGK) is used to generate the bases for the solution. Given $\bH, \bd$ and changing preconditioner $\bL_k$, the FGK iterative process can be described as follows. Let $\bu_1 = \bd/||\bd||$ and $\bv_1 = \bH^{\top}\bu_1/||\bH^{\top}\bu_1||$. Then at the $k-$th iteration, we generates vectors $\bp_k, \bv_k$ and $\bu_{k+1}$ such that
\begin{equation}
\label{relationship}
    \bH\bP_k = \bU_{k+1}\bM_k \quad {\text and} \quad \bH^{\top}\bU_{k+1} = \bV_{k+1}\bT_{k+1},
\end{equation}
where $\bP_k = \begin{bmatrix}\bL^{-1}_1\bv_1 & \dots & \bL^{-1}_k\bv_k\end{bmatrix} \in \mathbb{R}^{n \times k}, \bM_k \in \mathbb{R}^{(k+1)\times k}]$ is upper Hessenberg, $\bT_{k+1} \in \mathbb{R}^{(k+1)\times(k+1)}$ is upper triangular, and $\bU_{k+1} = \begin{bmatrix}\bu_1 & \dots & \bu_{k+1}\end{bmatrix} \in \mathbb{R}^{m \times (k+1)}$ and $\bV_{k+1} = \begin{bmatrix}\bv_1 &  \dots & \bv_{k+1}\end{bmatrix} \in \mathbb{R}^{n \times (k+1)}$ contain orthonormal columns. We remark that the column vectors of $\bP_k$ don't span a Krylov subspace like conventional Golub-Kahan bidiagonalization process\cite{o1981bidiagonalization, bjorck1988bidiagonalization}, but they do provide a basis for the solution $\bs_k$ at $k-$th iteration. Given the relationships in (\ref{relationship}), an approximate least-squares solution can be computed as $\bs_k = \bP_k\bq_k$, where $\bq_k$ is the solution to the projected least-squares problem,
\begin{equation}
\label{prob:project}
    \underset{\bs_k \in \mathcal{R}(\bP_k)}{\min}||\bH\bs_k - \bd||^2_2 = \underset{\bq_k }{\min}||\bM_k\bq_k - \beta_1\be_1||^2_2,
\end{equation}
where $\be_1 \in \mathbb{R}^{k+1}$ is the first column of $(k+1)$ by $(k+1)$ identity matrix. Although it is well known that iterative methods such as LSQR, it is well known that, for inverse problems, iterative methods without an appropriate regularization term exhibit semiconvergence behavior whereby the reconstructions eventually become contaminated
with noise and errors. Thus, a standard regularization term is included in (\ref{prob:project}), so that 
\begin{equation}
    \bq_k = \arg \underset{\bq}{\min}||\bM_k\bq - \beta_1\be_1||^2_2 + \lambda||\bq||^2_2.
\end{equation}
Henceforth, $\bs_k = \bP_k\bq_k$ is the numerical solution at $k-$th iteration for the full problem. To get a better regularized solution, we consider using weighted generalized crossed validation (WGCV) method \cite{chung2008weighted} to choose $\lambda$. 

%\subsection{Convergence}

%Here we state our main result as \cref{thm:bigthm}; the proof is
%deferred to \cref{sec:proof}.

%\section{Algorithm}
%\label{sec:alg}

%%\lipsum[40]

%Our analysis leads to the algorithm in \cref{alg:buildtree}.

%\begin{algorithm}
%\caption{Build tree}
%\label{alg:buildtree}
%\begin{algorithmic}
%\STATE{Define $P:=T:=\{ \{1\},\ldots,\{d\}$\}}
%\WHILE{$\#P > 1$}
%\STATE{Choose $C^\prime\in\mathcal{C}_p(P)$ with $C^\prime := \operatorname{argmin}_{C\in\mathcal{C}_p(P)} \varrho(C)$}
%\STATE{Find an optimal partition tree $T_{C^\prime}$ }
%\STATE{Update $P := (P{\setminus} C^\prime) \cup \{ \bigcup_{t\in C^\prime} t \}$}
%\STATE{Update $T := T \cup \{ \bigcup_{t\in\tau} t : \tau\in T_{C^\prime}{\setminus} \mathcal{L}(T_{C^\prime})\}$}
%\ENDWHILE
%\RETURN $T$
%\end{algorithmic}
%\end{algorithm}

%\lipsum[41]

\section{Numerical Results}
\label{sec:experiments}
In this section, we have two types of numerical experiments for testing the performance of our algorithm. In all the simulations, the initial guesses are randomly generated. The stopping criterion used in the all experiments depends on two parameters: one is the upper bound of the number of iterations $m_{\text{max}}$, and the other is the tolerance $\epsilon_{\text{tol}}$ of the relative difference between the observation and the approximation to decide whether the convergence has been achieved. The regularization parameter $\lambda$ is iteratively updated by the flexible Krylov method with weight generalized cross validation method. These experiments ran on a laptop computer with Intel i5 CPU 2GHz and 16G memory.

\subsection{Image recovering by tensor completion}
For the first experiment, we test two cases for this example, where the missing pixels for the first case are randomly chosen while the miss part for the second case is deterministic. The reconstruction error is computed with the relative error $\frac{||\mathcal{A}_{\text{estimate}} - \mathcal{A}||_2}{|| \mathcal{A} ||_2}$, where $\mathcal{A}_{\text{estimate}}$ denotes the approximated tensor, and $\mathcal{A}$ represents the tensor we want to reconstruct. 

\begin{figure}[h!]
    \begin{tabular}{cc}
    \hspace{-0.5in}
      \includegraphics[width = 0.56\textwidth]{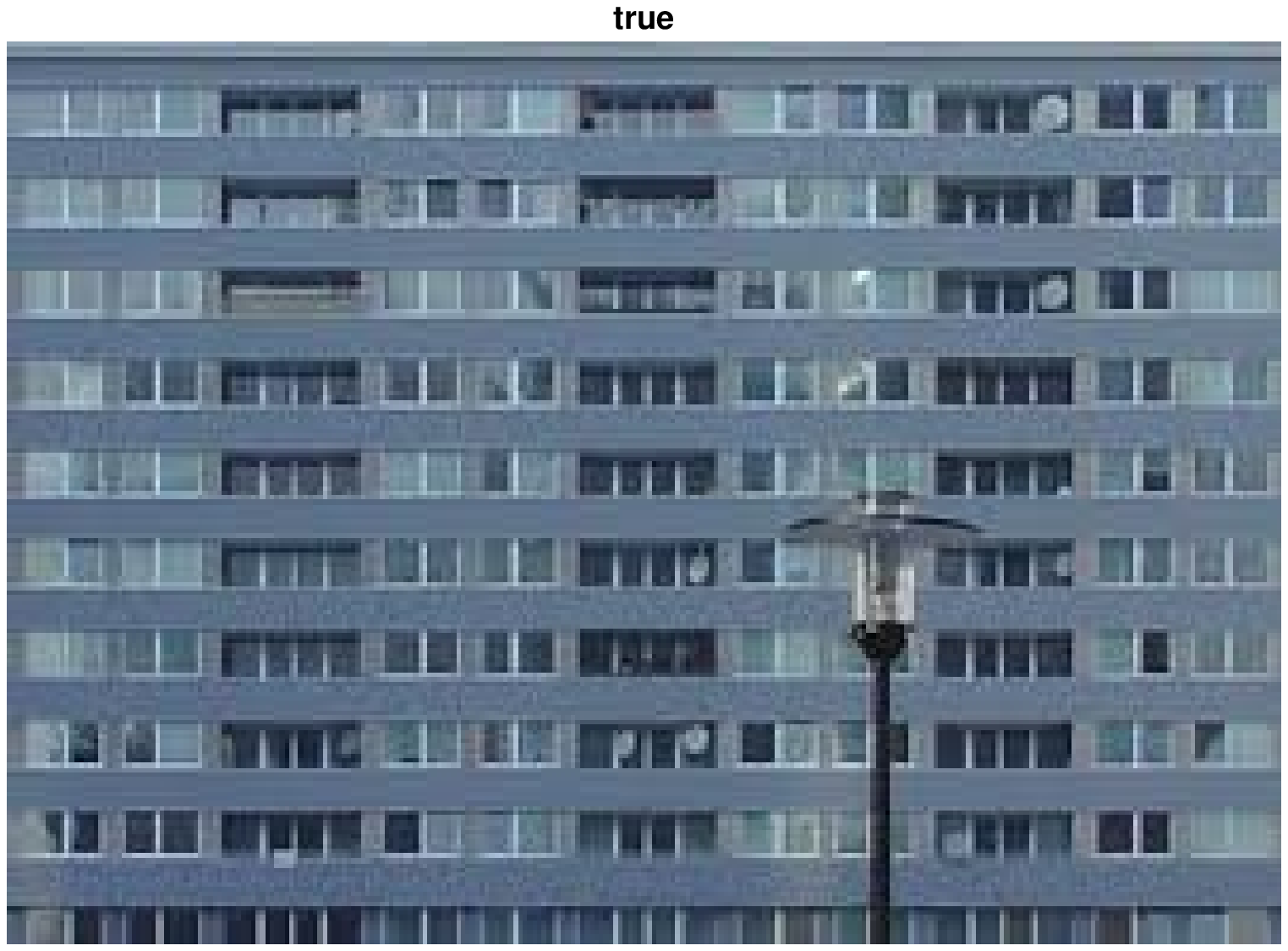}& \hspace{-0.5in}\includegraphics[width = 0.56\textwidth]{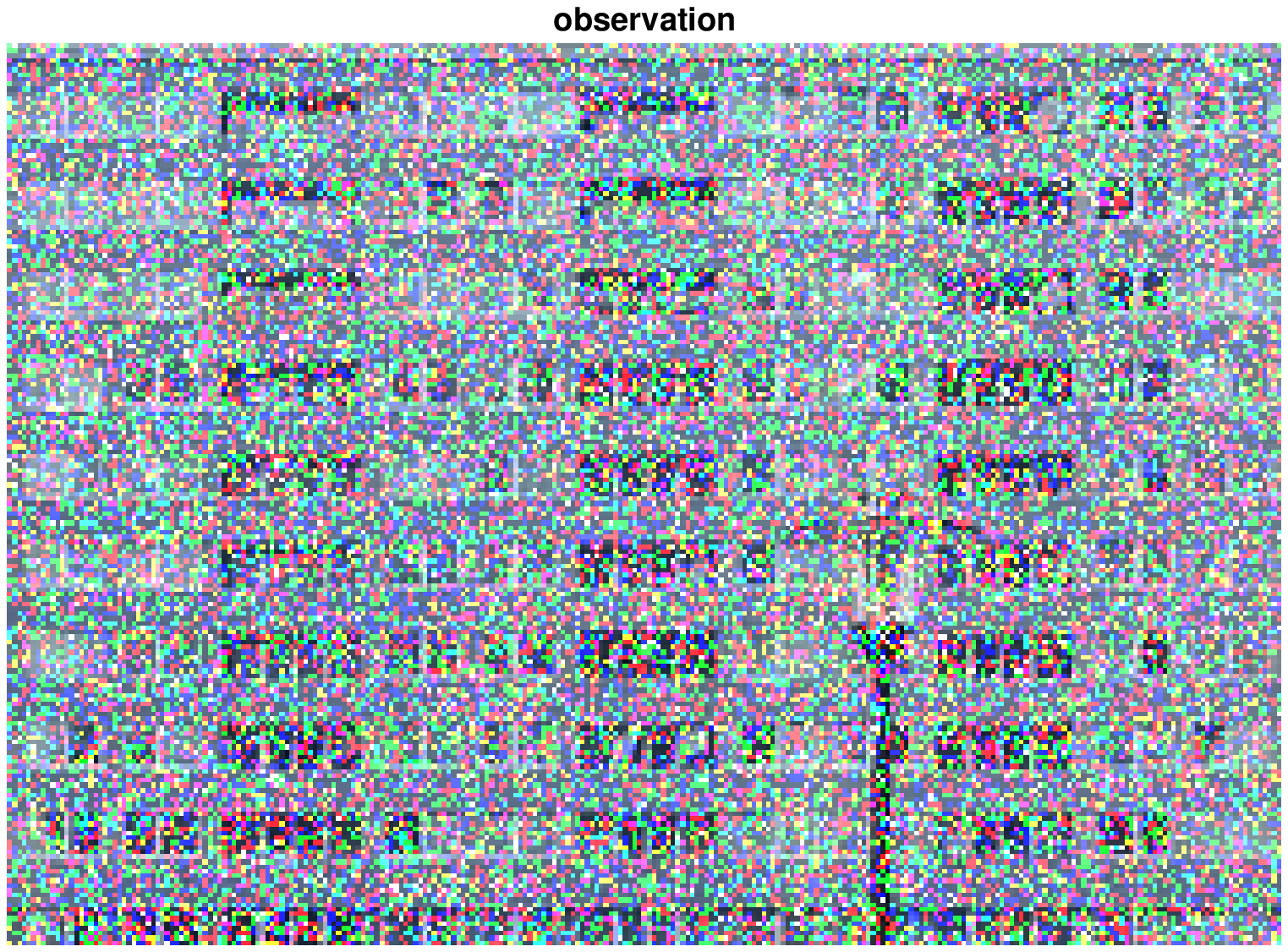}
    \end{tabular}
    \caption{True image (left) and observation (right). }
    \label{fig:trueob}
\end{figure}

{\bf Case 1:} We implemented our algorithm on a color image $\mathcal{A} \in \mathbb{R}^{189 \times 267 \times 3}$ shown in Figure \ref{fig:trueob}. We recovered an estimated color image after removing $30\%$ of the entries from the origin color image, which is shown in Figure \ref{fig:trueob}. The upper bound $R$ of rank$(\mathcal{A})$ is fixed to 50 in the algorithm. The stopping criteria for this case are assumed to be $m_{\text{max}} = 500$ and $\epsilon_{\text{tol}} = 10^{-3}$. We choose $\lambda = 35$ for the conventional CP tensor. 

\begin{figure}[h!]
    \begin{tabular}{cc}
    \hspace{-0.5in}
      \includegraphics[width = 0.56\textwidth]{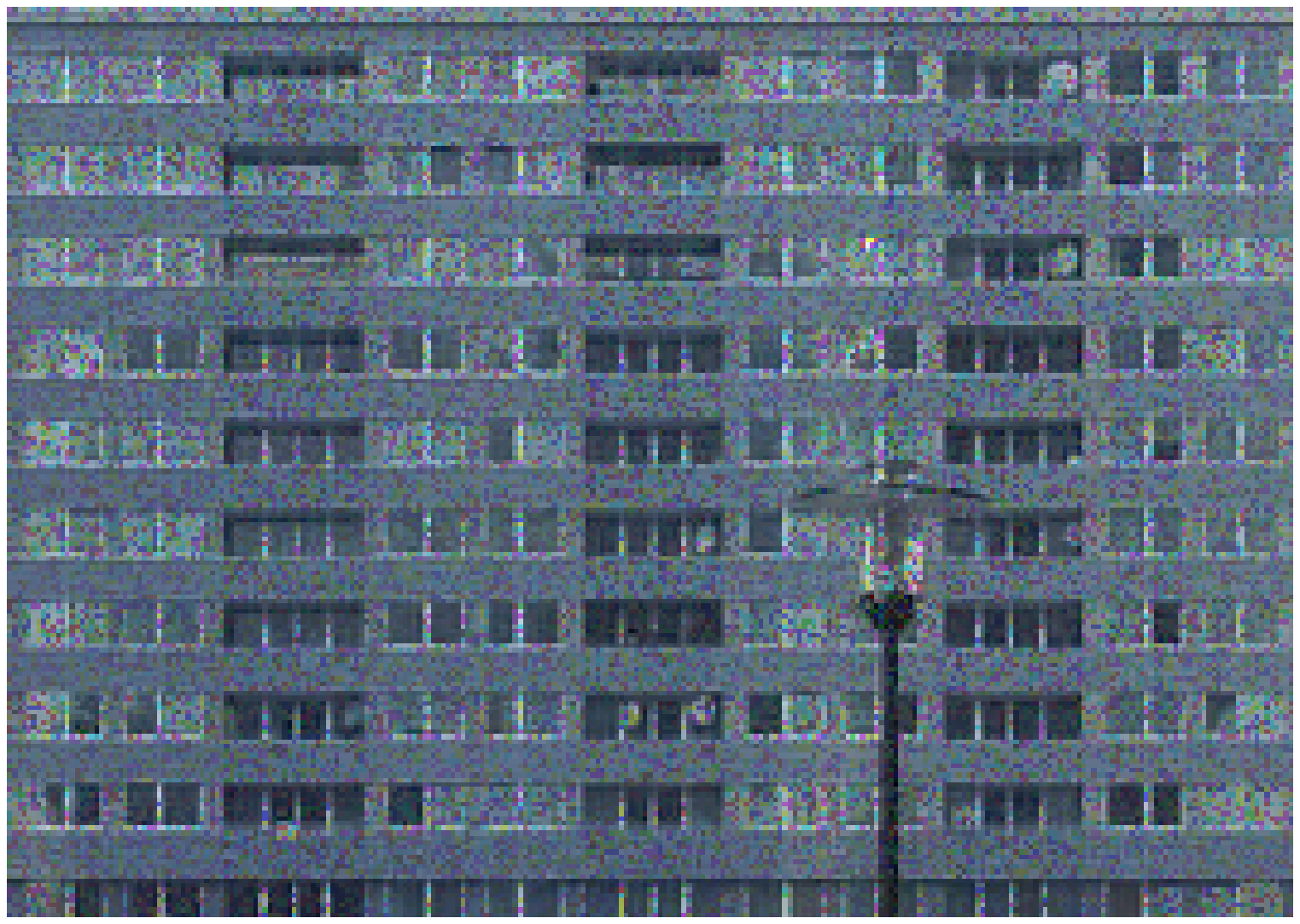}& \hspace{-0.5in}\includegraphics[width = 0.56\textwidth]{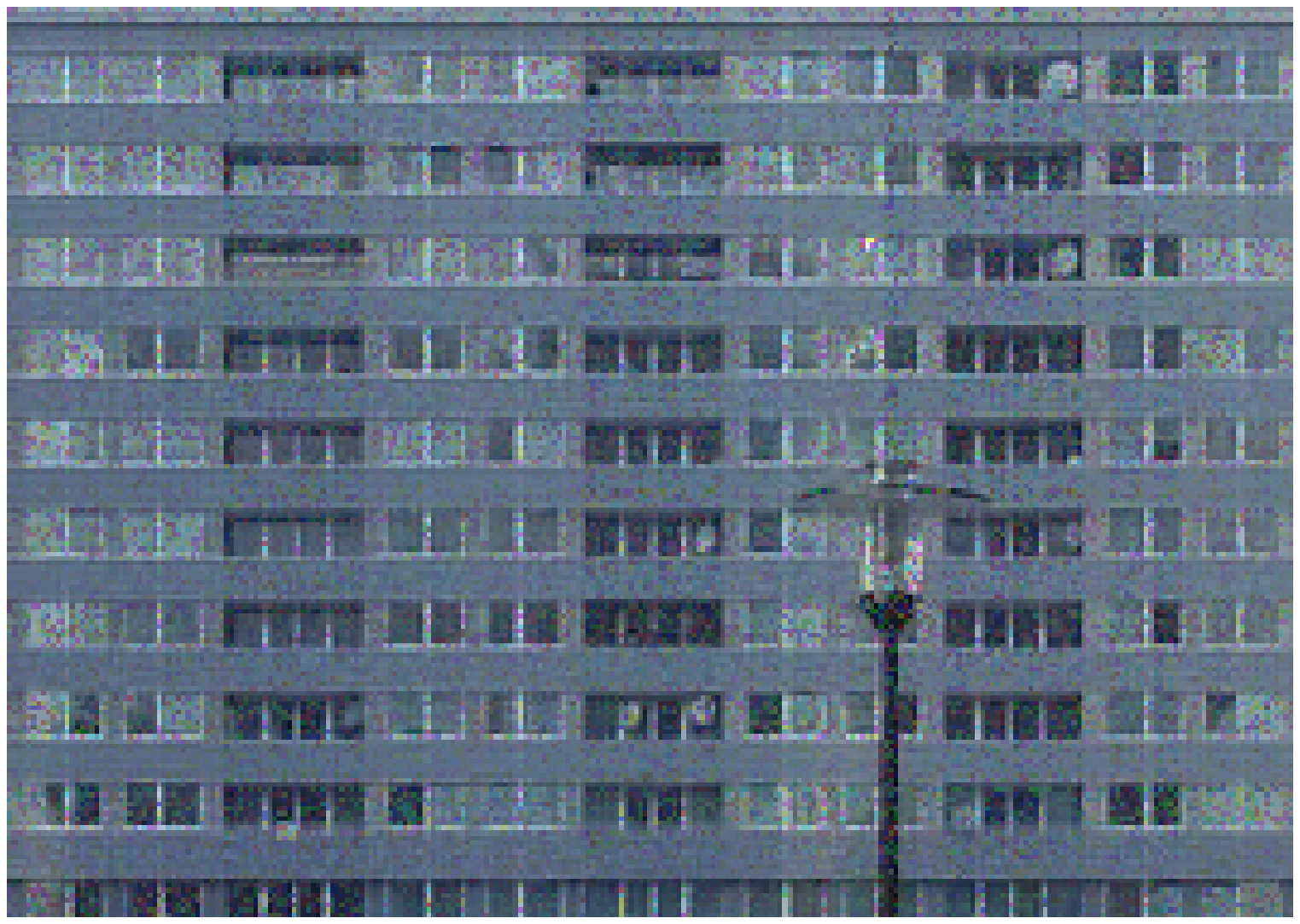}\\
       CP tensor(0.1921)& \hspace{-0.5in} practical regularization CP tensor(0.1296)
    \end{tabular}
    \caption{Recovered color image with relative reconstruction error norms provided in the bracket. }
    \label{tab: reconstruction}
\end{figure}

The recovered images by original CP tensor and practical regularization CP tensor are provided in Figure \ref{tab: reconstruction}. We can see that practical regularization CP tensor produces recovered image that has much less noise than classical CP tensor, demonstrating that using flexible Krylov method to determine different regularization parameter for each iteration is beneficial.  The comparison of the relative error shown in Figure \ref{tab: reconstruction} also verifies the better performance of our practical regularization CP tensor. 

\begin{figure}[h!]
    \begin{tabular}{cc}
      \includegraphics[width = 0.55\textwidth]{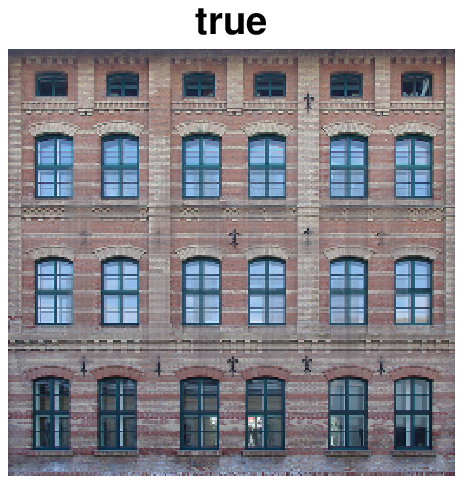}&
      \hspace{-0.9in}
      \includegraphics[width = 0.55\textwidth]{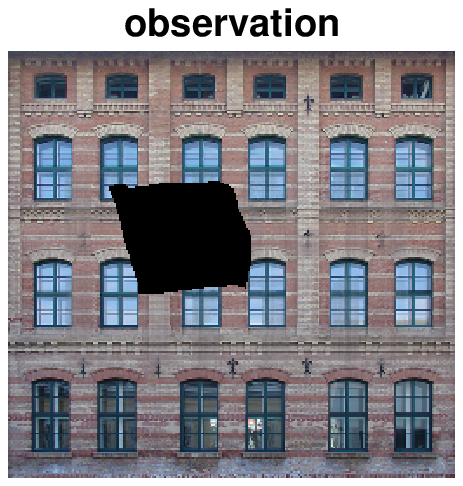}
    \end{tabular}
    \caption{True image (left) and observation (right).}
    \label{fig:trueob_uab}
\end{figure}

{\bf Case 2: }We consider recovering the image $\mathcal{A} \in \mathbb{R}^{246\times 257\times 3}$ with the certain missing pixels as shown in Figure \ref{fig:trueob_uab}, associating with its true image. The upper bound $R$ of rank$(\mathcal{A})$ is chosen to be $50$. The stopping criteria are setup as $m_{\text{max}} = 250$ and $\epsilon_{\text{tol}} = 10^{-3}$. For the classical CP tensor, we choose $\lambda  = 35$. 

\begin{figure}[h!]
    \begin{tabular}{cc}
      \includegraphics[width = 0.55\textwidth]{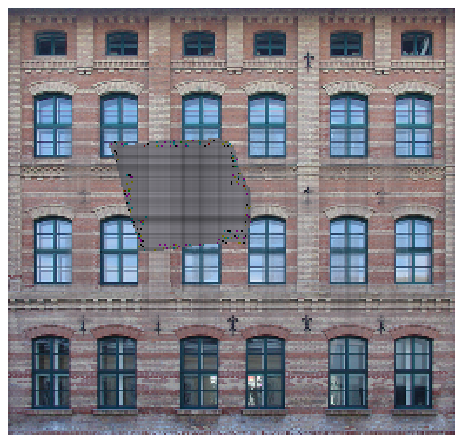}&
      \hspace{-0.9in}
      \includegraphics[width = 0.55\textwidth]{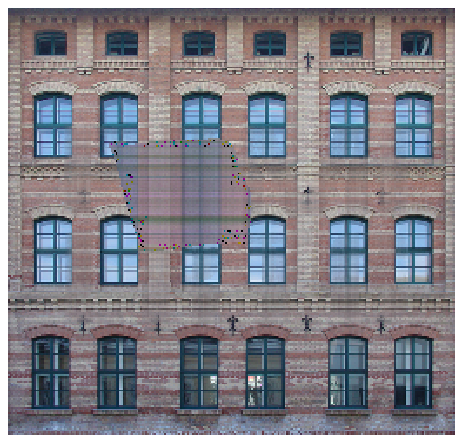}\\
       CP tensor(0.0950)& \hspace{-0.7in} practical regularization CP tensor(0.0743)
    \end{tabular}
    \caption{Recovered color image with relative reconstruction error norms provided in the bracket. }
    \label{tab: reconstruction_uab}
\end{figure}

The recovered images are provided in Figure \ref{tab: reconstruction_uab}. We observe that our algorithm does obviously more complete recovering than conventional CP tensor, which is also demonstrated by the relative error.

\subsection{Model order reduction}
Next we investigate a scenario in model order reduction where the key snapshots needs to be obtained to capture the low rank structure of the solution manifold that has low Kolmogorov width \cite{lorentz1996constructive, pinkus2012n}. This example demonstrates advantages of our practical regularization CP tensor. 
%We consider two different cases. 
%\begin{enumerate}
    %\item 
    Model order reduction techniques such as the POD and the Reduced Basis Methods are typically used to solve the problems requiring one to query an expensive yet deterministic computational solver once for each parameter node. We shows that hybridizing our approach and regularized alternating least-squares method \cite{navasca2008swamp, li2013some} provides a novel way to to do model reduction and pattern extraction. More specifically, assuming $\mathcal{A}$ is the collection of the solutions on sampled parameters. To select the snapshots (reduced bases) for the low rank approximation of the solution manifold, we employ our algorithm to give a prior knowledge of rank($\mathcal{A}$) denoted as $R$, and then we run regularized alternating least-squares method according to the $R$ to approximate $\mathcal{A}$ and build up the reduced bases.

    % \item Many computational methods for stochastic problems in uncertainty quantification (UQ) are accelerated by model order reduction algorithms, and much recent effort in this direction has been rewarded with many promising developments [CITE]. Motivated by the tensor recovery based model reduction techniques proposed in \cite{Tang2020RankAT}, we implement our approach in tensor recovery to approximate the generalized Polynomial Chaos (gPC) [CITE] coefficients. 
%\end{enumerate}

%\subsubsection{Tensor Approximation}
In this experiment, we consider the following two-dimensional diffusion equation n that induces a solution manifold that requires
many more snapshots to achieve small error:
\begin{equation}
    (1+\mu_1x)u_{xx} + (1+\mu_2y)u_{yy} = e^{4xy} \qquad \text{on} \quad \Omega.
\end{equation}
The physical domain is $\Omega = [-1, 1]\times [-1, 1]$ and we impose homogeneous Dirichlet boundary conditions on $\partial \Omega$. The truth approximation is a spectral Chebyshev collocation method \cite{hesthaven2007spectral, trefethen2000spectral} with $\mathcal{N}_x = 100$ degrees of freedom in each direction.  This means the truth approximation has dimension $\mathcal{N} = 10000 (\mathcal{N} = \mathcal{N}^2_x)$. 
The parameter domain $\mathcal{D}$ for $(\mu_1, \mu_2)$ is taken to be $[-0.99, 0.99]\times [-0.99, 0.99]$. For the parameter sample set, we discretize $\mathcal{D}$ using a tensorial $9 \times 9$ cartisian grid, thus the size of training is $81$. The testing set $\Xi_{\text{test}}$ contains another $10$ random samples in $\mathcal{D}$.
The resulting tensor $\mathcal{A}$ is of dimension $100 \times 100 \times 81$. Given an initial value of the rank $R_0= 50$ and tolerance $\varepsilon= 10^{-2}$, we run our algorithm on $\mathcal{A}$ at first, then sort the rescaling coefficients $\{\alpha_r\}_{r = 1}^{R_0}$ in descending order and discard the the coefficients below than $\varepsilon\alpha_{\text max}$, where $\alpha_{\text max}$ is the maximum  rescaling coefficients. Assuming that the number of the coefficients we keep is $R$ that is much smaller than $\mathcal{N}$, we run regularized alternating least-squares method with rank $R$ to approximate
\begin{equation}
    \mathcal{A} \approx \sum_{r=1}^{R} \alpha_{r}\mathbf{x}_r\circ\mathbf{y}_r\circ\mathbf{z}_r,
\end{equation}
and build up the reduced bases $\{\phi_r\}_{r=1}^{R}$ by orthonormalizing $\{\widehat{\phi}_r\}_{r=1}^{R}$, where $\widehat{\phi}_r$ is created by vectorizing 
$\mathbf{x}_r\circ\mathbf{y}_r$. 

\subsubsection{Computational performance}
Given initial value of $\lambda = 10$, our proposed approach chooses $R = 20$ reduced bases, while the classical CP tensor chooses $R = 7$ reduced bases with the same initialization. In addition,  in this model reduction example, the number of reduced bases built up by the classical CP tensor is very sensitive to the choice of initial value of $\lambda$, that is, classical CP tensor will chooses much less(more) reduced bases if we increase(decrease) $\lambda$ a little bit. This means that our practical regularization CP tensor including the automatic regularization parameter selection is more robust than the conventional CP tensor, since it doesn't rely on the empirical choice of the initial value of the regularization parameter. 

\begin{figure}[h!]
   \centering
       \begin{tabular}{cc}
       \includegraphics[width = 0.48\textwidth]{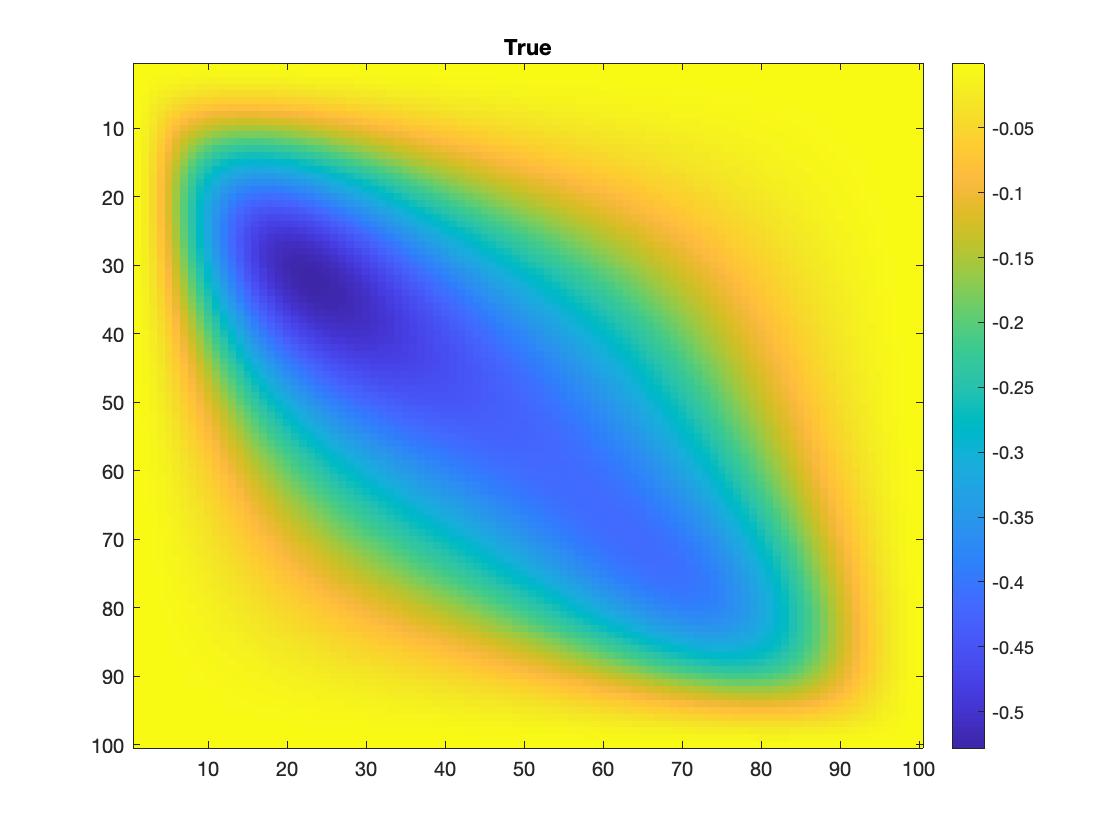}  & 
    %   \hspace{-0.3in} \includegraphics[width = 0.33\textwidth]{true_2.jpg} & 
       \hspace{-0.3in} \includegraphics[width = 0.48\textwidth]{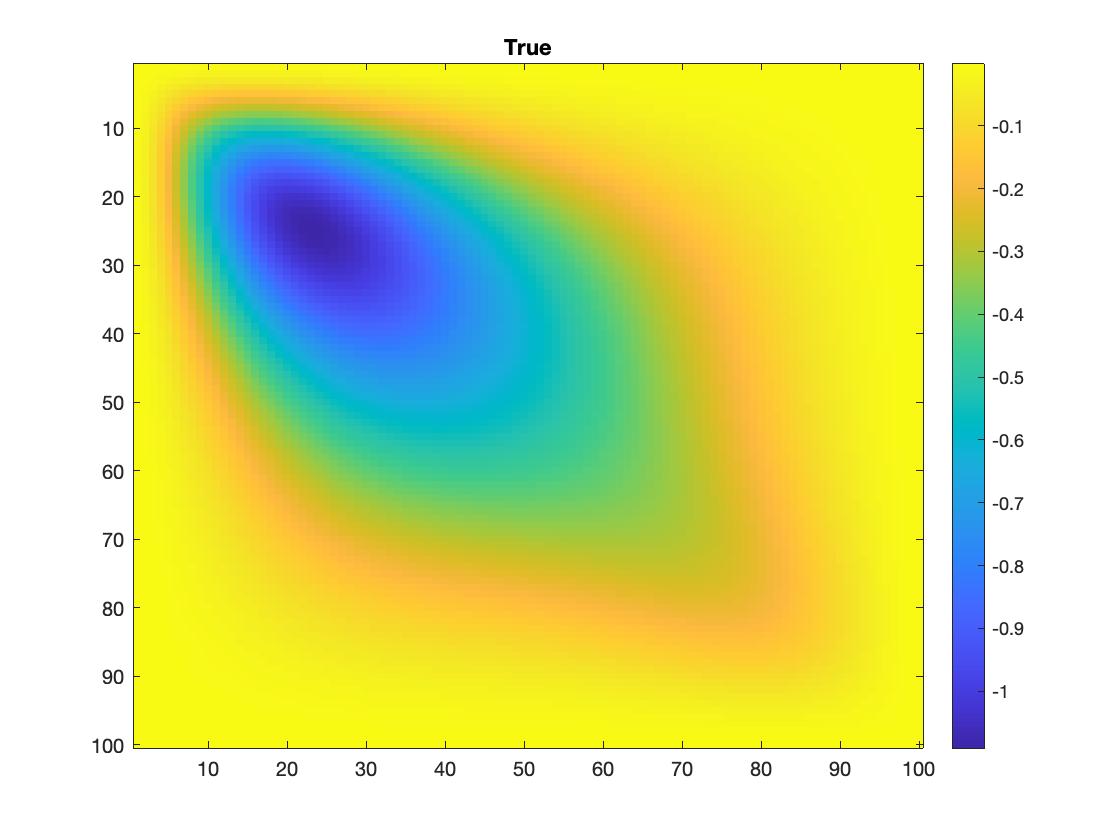}\\
           \includegraphics[width = 0.48\textwidth]{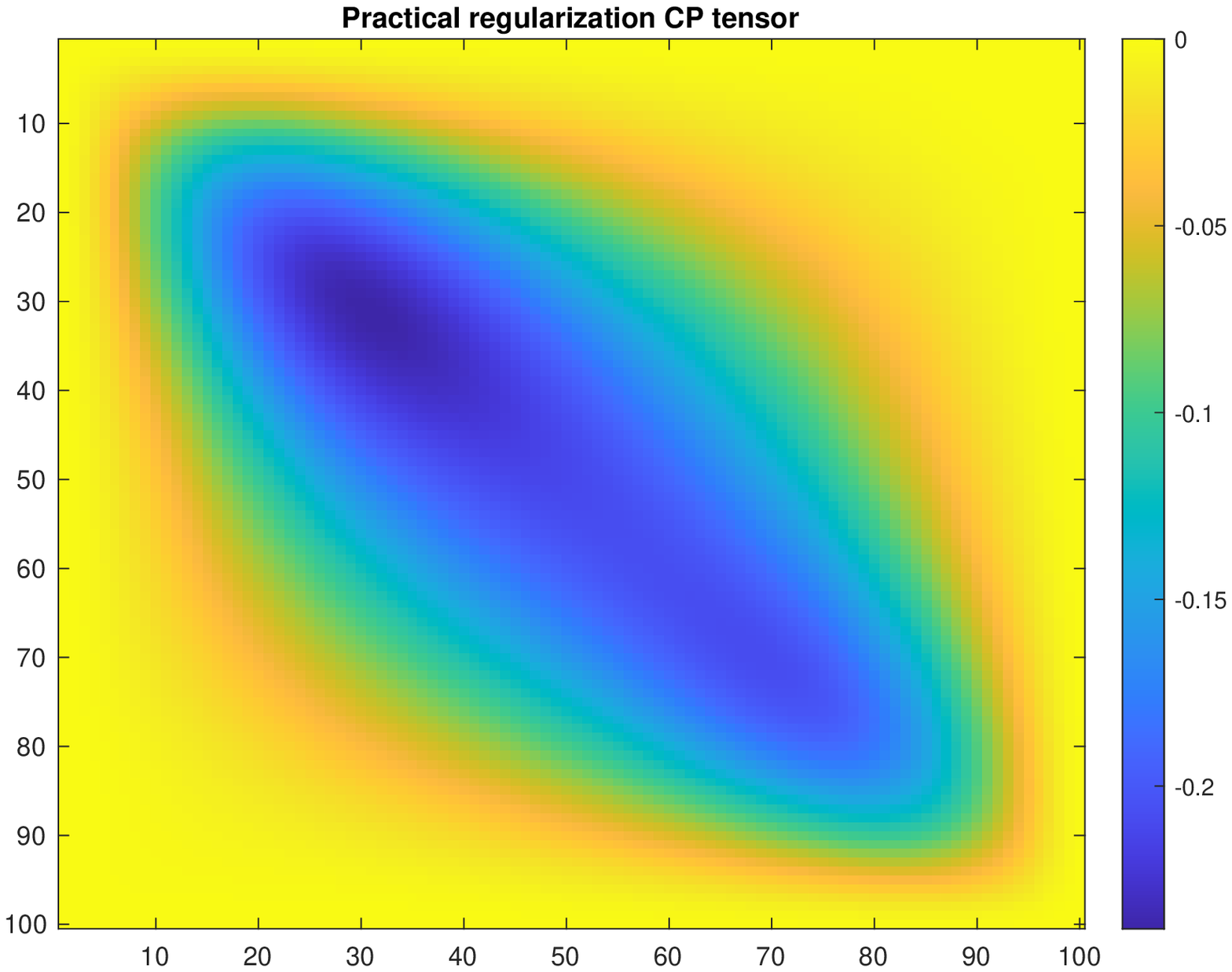}  
           %&  \hspace{-0.3in} \includegraphics[width = 0.33\textwidth]{tensor_2.jpg}
        &  \hspace{-0.3in} \includegraphics[width = 0.48\textwidth]{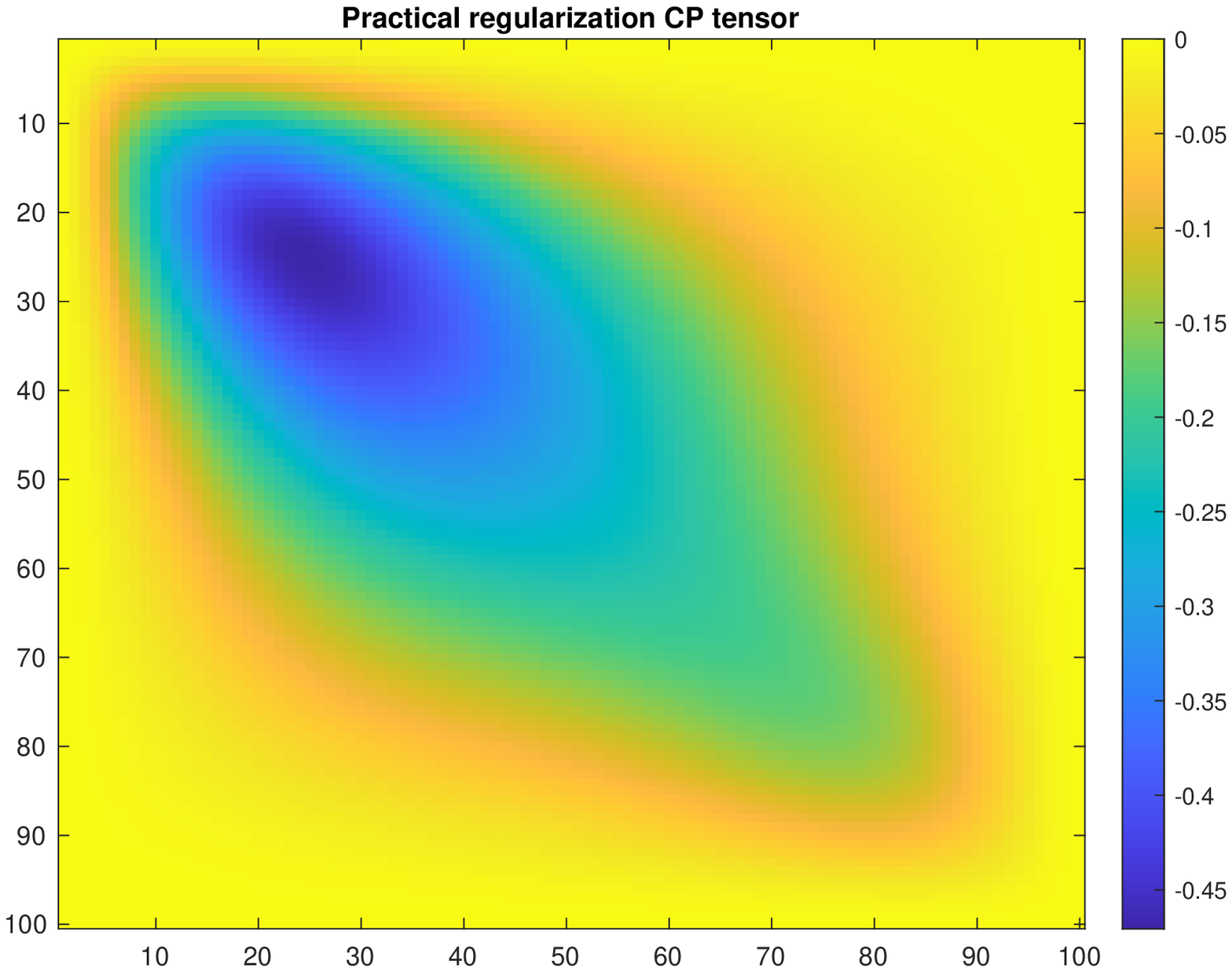}\\
          $\ell^2$ error: $3.4\times 10^{-1}$
          %& $L^2$ error: $3.6\times 10^{-1}$ 
          & $\ell^2$ error: $5.9\times 10^{-1}$\\
         \includegraphics[width = 0.48\textwidth]{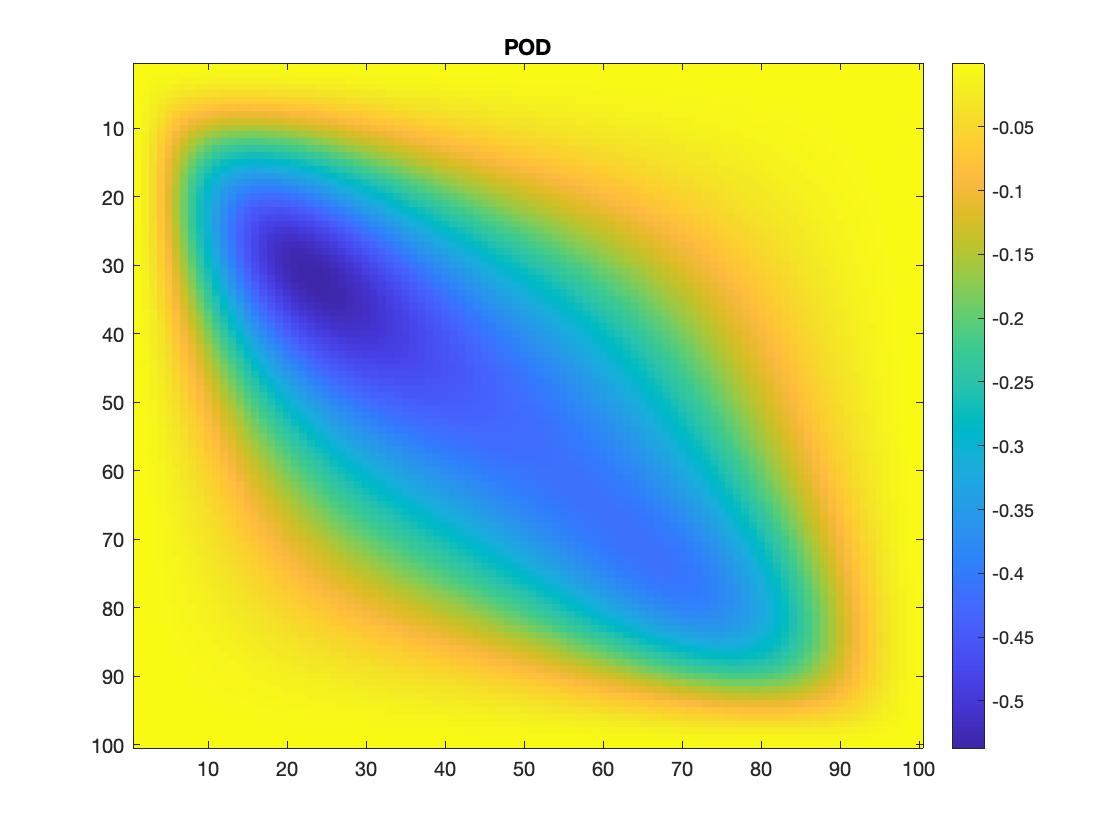}  &  \hspace{-0.3in} %\includegraphics[width = 0.33\textwidth]{POD_2.jpg}&  \hspace{-0.3in} 
         \includegraphics[width = 0.48\textwidth]{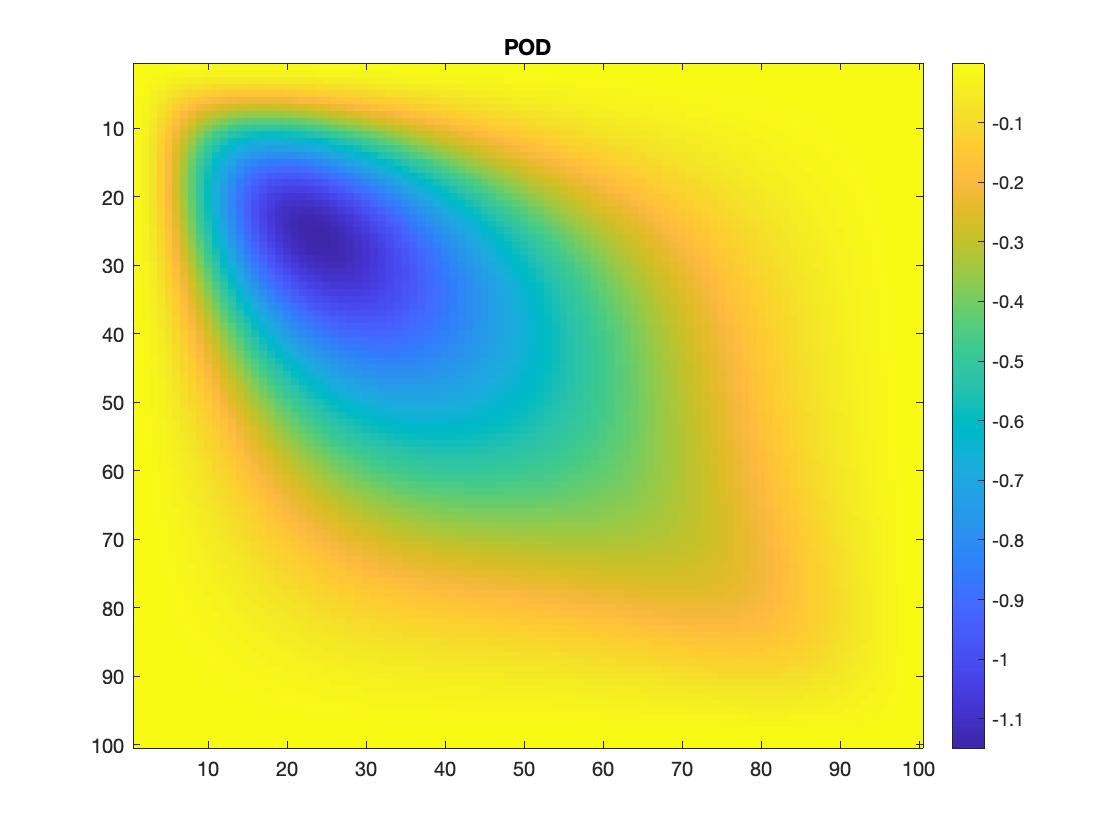}
    \\
        $\ell^2$ error: $1.3\times 10^{-2}$
        %& $L^2$ error: $3.0\times10^{-2}$
        & $\ell^2$ error: $6.8\times10^{-2}$\\
       \end{tabular}
       \caption{The approximation of the solution of Two-dimensional diffusion equation at two parameters in the testing set: $\mu_1 = -0.9193, \mu_2 = 0.6913$(left) and $\mu_1 = -0.8986,   \mu_2 = -0.7977$(right)}
       \label{fig:recon-mor}
   \end{figure}

\begin{figure}[h!]
       \centering
       \includegraphics[width = \textwidth]{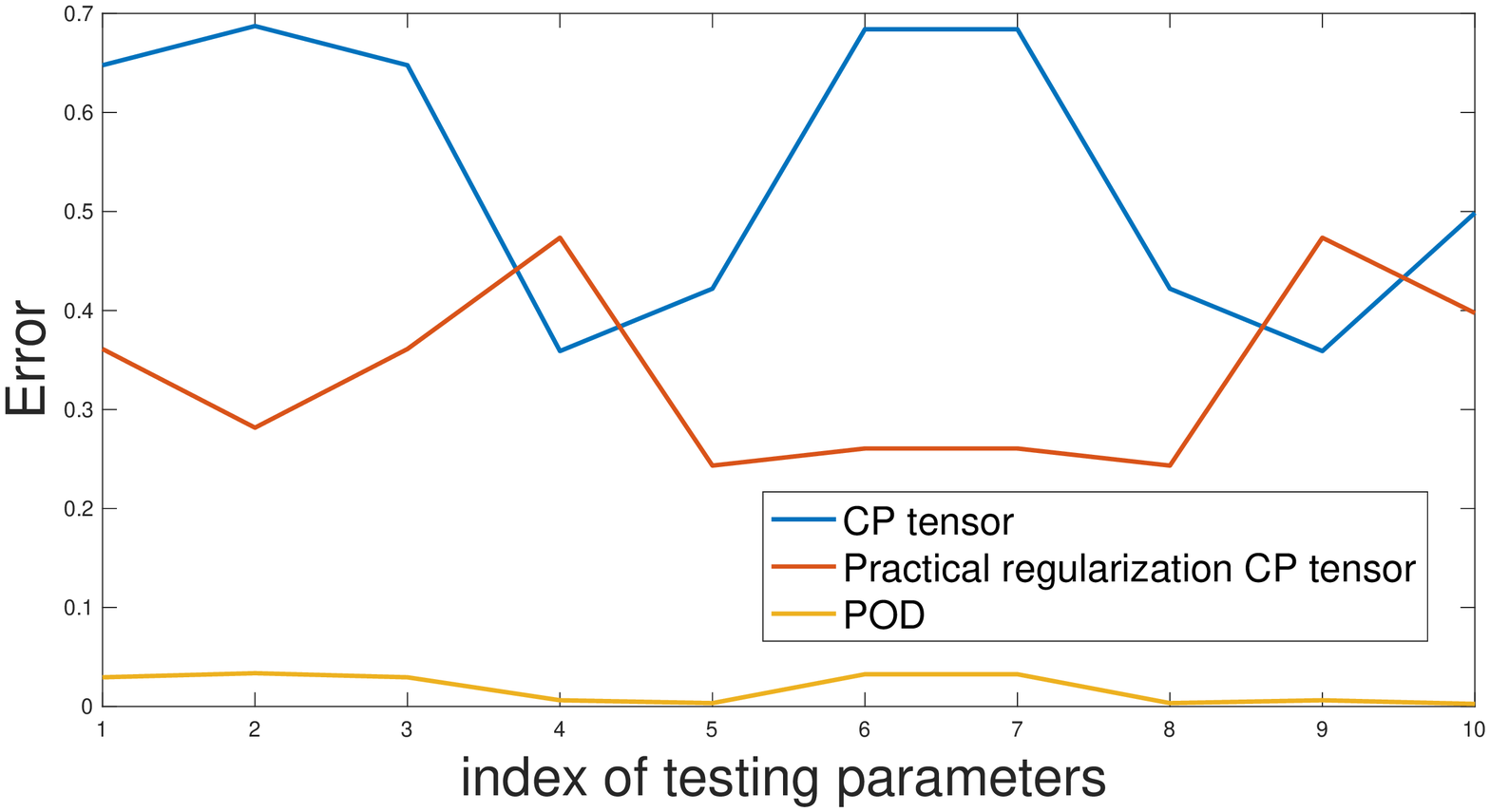}
       \caption{The approximation of $u(x, \mu_1, \mu_2)$ with $(\mu_1, \mu_2) \in \Xi_{\text{test}}$: CP tensor, practical regularization CP tensor and POD.}
       \label{fig:error-mor}
   \end{figure}

To gain the understanding of the quality of the reduced bases our proposed approach, the number of reduced bases constructed by the POD is also set to be $20$ for comparative purposes.  The algorithmic accuracy is evaluated by approximating the $u(x, \mu_1, \mu_2)$ via reduced bases, where $(\mu_1, \mu_2) \in \Xi_{\text{test}}$ and the error is measured by $\ell^2$ norm. Figure \ref{fig:recon-mor} displays the approximation of the solution at two parameters drawn from the testing set for the practical regularization CP tensor and the POD. We observe that our proposed algorithm faithfully captures the feature of the solution, although the error is large compared to the POD. The performance of approximating the solutions at all the parameters in the testing set for our algorithm, canonical CP tensor and the POD is provided in Figure \ref{fig:error-mor}. Table \ref{tab:mor} further quantifies the range of the approximation results. It can clearly be seen that approximation quality of practical regularization CP tensor is better than original CP tensor due to that the number of bases adopted by practical regularization CP tensor doesn't depend on the initial value of the regularization parameter. However, the approximation error of our algorithm is large compared with the error resulting from the POD. 

 \begin{table}[h!]
   \hspace{-0.2in}
       \begin{tabular}{cccc}
       \toprule
        & CP tensor  & practical regularization CP tensor  &  POD \\ 
        \midrule
       Error &$[0.36, \; 0.83]$ & [0.24, \; 0.58] & [0.003, \; 0.068]\\ 
        \midrule
        Number of bases & $7$ & $20$ & $20$ \\
        \bottomrule
       \end{tabular}
       \caption{Comparison of algorithmic  accuracy of approxmating $u(x, \mu_1, \mu_2)$ with $(\mu_1, \mu_2) \in \Xi_{\text{test}}$ for CP tensor, practical regularization, CP tensor and POD}
       \label{tab:mor}
   \end{table}
   
\subsubsection{A note on compression}
For this model reduction problem, although both practical regularization CP tensor and CP tensor provide a more parsimonious representation of the data than the POD, comparing the compression ratios between the CP tensor techniques and the POD illustrates the difference and the benefit of the  CP tensor techniques. For a rank $R = 20$ tensor of dimension $\mathcal{A} = 100 \times 100 \times 81$, the compression, the compression ratios are
\begin{align}
    \mathcal{C}_{\text{POD}} &= \frac{I \cdot J \cdot K}{R\cdot(I\cdot J + K + 1)} = \frac{100^2\cdot 81}{20\cdot(100^2 + 81 + 1)} \approx 4.02, \\
    \mathcal{C}_{\text{CP}} & = \frac{I \cdot J \cdot K}{R\cdot(I\cdot J + K + 1)} = \frac{100^2\cdot 81}{20\cdot(100+100 + 81 + 1)} \approx 143.62. 
\end{align}
Notice that the POD requires the tensor to be reshaped in some direction. The comparison illustrates the striking difference between the compression ratios. It is worth mentioning that the CP tensor approaches requires much less memory to approximate the data. This can be of importance if the online stage (approximating the data) is in limited storage situations and that the accuracy requirement is not high. 

%\lipsum[50]

% \Cref{fig:testfig} shows some example results. Additional results are
% available in the supplement in \cref{tab:foo}.

% \begin{figure}[htbp]
%   \centering
%   \label{fig:a}\includegraphics{lexample_fig1}
%   \caption{Example figure using external image files.}
%   \label{fig:testfig}
% \end{figure}

%\lipsum[51]

% \section{Discussion of \texorpdfstring{{\boldmath$Z=X \cup Y$}}{Z = X union Y}}

% %\lipsum[76]

% \section{Two}

% \section{Three}

\section{Conclusion}
\label{sec:conclusions}
In this paper, we have presented a new low-rank CP tensor completion algorithm by combining the flexible hybrid method and the CP tensor completion. A key advantage of this method is that the regularization parameter can be easily and automatically estimated during the iterative process, which substantially reduces the difficulty of  initializing the regularization parameter and improves the robustness of the algorithm.  In addition to memory savings, our proposed approach demonstrates outstanding performance on the model reduction example, compared to the POD. Moreover, our image recovery experiments show that our algorithm has a practical advantage in capturing more details in image reconstruction over the conventional CP tensor due to a more optimal choice of the regularization parameter. In our future outlook, we will extend this hybrid approach in a tensor based total variation formulation for denoising and deblurring multi-channel images and videos.

\appendix
\section{Proximal Gradient} 
%\lipsum[71]

Recall $g(\sigma)= \lambda \parallel \sigma\parallel_{1}$ is convex and non-differentiable. The function
g can be turn into a proximal operator to find its minimum using the definition \cite{ParikhBoyd,7938377} below:
\begin{definition}
Given a proper closed convex function f: $\mathbb{R}^{n}\longrightarrow \mathbb{R}\bigcup {\infty}$, the proximal operator scaled by $\delta>0$, is a mapping from $\mathbb{R}^{n}\longrightarrow \mathbb{R}$ defined by 
\begin{equation*} 
\textbf{prox}_{\delta g}(v):=
\begin{aligned}
& \underset{y \in \mathbb{R}^{n} }{\text{argmin}}
& &( g(y) + \frac{1}{2\delta} \parallel y-v\parallel^{2} ).
\end{aligned}
\end{equation*}
\end{definition}

\noindent Then the proximal operator for $g(\sigma)$ is,
~\vspace*{-8pt}
\begin{equation*}
\begin{split}
prox_{g}(v)
& = \underset{\sigma}{\text{argmin }}  (g(\sigma) + \frac{1}{2\delta}\parallel\sigma - v\parallel_{2}^{2})\\
& =\underset{\sigma}{\text{argmin }} (\lambda\parallel\sigma\parallel_{1} + \frac{1}{2\delta} \parallel\sigma - v\parallel_{2}^{2})\\
& = \underset{\sigma}{\text{argmin }} (\lambda \sum_{i=1}^{n} \mid\sigma_{i}\mid + \frac{1}{2\delta}\sum_{i=1}^{n} (\sigma_{i} - v)^{2})\\
\end{split}
\end{equation*}
~\vspace*{-12pt}
For $\sigma = (\sigma_1 ,..., \sigma_n)$
\begin{multline*}
(prox_{\delta g}(v)_{i})_{i}
 = \underset{\sigma_i}{\text{argmin }} (\lambda\mid\sigma_i\mid + \frac{1}{2\delta}(\sigma_i - v_i)^{2})\\
 = \underset{\sigma_i}{\text{argmin }} \begin{cases}
(\lambda - \frac{v_i}{\delta})\sigma_i + \frac{1}{2\delta}\sigma_i^{2} , \sigma_i>0
\\
-(\lambda + \frac{v_i}{\delta})\sigma_i + \frac{1}{2\delta}\sigma_i^{2},  \sigma_i<0 \Bigg\}
\end{cases}
\end{multline*}
%Case 1: for $\lambda - \frac{v_i}{\delta}<0$, $\Rightarrow v_i>\lambda\delta$~\\
%\text {The minimum is when }$\sigma_i>0$ then the derivative will be $\lambda -\frac{v_i}{\delta} + \frac{1}{\delta}\sigma_i =0 \Rightarrow \sigma_i = v_i - \lambda\delta$~\\
%Case 2: for $\lambda + \frac{v_i}{\delta}<0 \Rightarrow v_i<-\lambda\delta$~\\
%\text{Then the minimum occurs when } $\sigma_i<0$ \text{ and }$\sigma_i = v_i + \lambda\delta$~\\
%Case 3: for $-\lambda\leqslant \frac{v_i}{\delta}\leqslant\delta$, \text{ the minimum should be 0}.~\\
Hence, \begin{equation*}
prox_{g}(v)_{i} =
\begin{cases}
0 ,                           &\mid v_{i}\mid\leq\lambda\\
v_{i} - \lambda sign(v_{i}),  &\mid v_{i}\mid>\lambda
\end{cases}
\end{equation*} with $\delta=1$\\

\section*{Acknowledgments}
This material is based upon work supported by the National Science Foundation under Grant No. DMS-1439786 while the authors, J. Jiang and C. Navasca, were in residence at the Institute for Computational and Experimental Research in Mathematics in Providence, RI, during the Model and Dimension Reduction in Uncertain and Dynamic Systems Program.

\bibliographystyle{siamplain}
\bibliography{tensor_reg}
\end{document}